%% file: prime.tex
\documentclass{amsart}
\usepackage{amssymb, epsfig}

\newtheorem{thm}{Theorem}[section]
\newtheorem{cor}[thm]{Corollary}
\newtheorem{lem}[thm]{Lemma}
\newtheorem{prop}[thm]{Proposition}

\newtheorem{remarkk}[thm]{Remark}
\newenvironment{remark}{\begin{remarkk} \em}{\end{remarkk}}
\newtheorem{examplee}[thm]{Example}

\newtheorem{defnn}[thm]{Definition}

\newcommand{\bbz}{\mathbb{Z}}
\newcommand{\bbr}{\mathbb{R}}

\title{Invertible Knot Concordances and Prime Knots}
\author{Se-Goo Kim}
\address{Department of Mathematics, Indiana University, Bloomington,
Indiana 47405}
\email{sekim@indiana.edu}
\date{March 5, 2000}


\begin{document}

\maketitle

\vspace*{-1em}

\section{Introduction}
Kirby and Lickorish \cite{kl79} showed that every knot in $S^3$ is
concordant to a prime knot, equivalently, every concordance class
contains a prime knot. Generalizations appear in
\cite{liv81,mye83,mye93,som84}.
Sumners \cite{sum71} introduced the notion of invertible concordance.
We prove here that the Kirby and  Lickorish's result can be strengthened:
\begin{thm}
\label{thm:main}
Every knot in $S^3$ is invertibly concordant to a prime knot.
\end{thm}
Corresponding to
invertible concordance there is a group, the {\sl double concordance
group}, studied in \cite{lev83,rub83,sto78}.
A consequence of our work is that every double concordance class contains
a prime knot.

\section{Definitions and basic results}
\label{sec:definitions}
In all that follows manifolds and maps will be smooth and orientable.
Let~$I$ denote the interval $[0,1]$.

A \emph{link} of $n$ components, $L$, is a smooth pair $(S^3,l)$
where $l$ is a smooth oriented submanifold of $S^3$
diffeomorphic to $n$ disjoint copies of $S^1$. 
A \emph{knot} $K$ is a link of one component.
Two links, $L_1$ and $L_2$, each of $n$ components, are called
\emph{concordant} if there exists a proper smooth oriented submanifold $w$
of $S^3\times I$, with $\partial w=(l_1\times 0\cup (-l_2)\times 1)$ and
$w$ diffeomorphic to $n$ disjoint copies of $S^1\times I$.
Let $(W;L_1,L_2)$ denote $(S^3\times I, w)$ the concordance between $L_1$
and $L_2$. If $(W_1;L_1,L_2)$ and $(W_2;L_2,L_3)$ are two concordances with
a common boundary component (oriented oppositely) we can then paste $W_2$
to $W_1$ along $L_2$ to get $(W_1\cup W_2;L_1,L_3)$.

A concordance $(W;L_1,L_2)$ is said to be \emph{invertible at} $L_2$ if
there is a concordance $(W';L_2,L_1)$ such that
$(W\cup W';L_1,L_1)$ is diffeomorphic to $(L_1\times I; L_1,L_1)$, the
product concordance of $L_1$.  Given the above situation, we say that
$L_1$ \emph{is invertibly concordant to} $L_2$,
and $L_2$ \emph{splits} $L_1\times I$.
In the same manner, concordance and invertible concordance
can be defined for knots and links in the solid torus $S^1\times D^2$.

A submanifold $N$ with boundary is said to be \emph{proper} in a
manifold $M$ if $\partial N=N\cap \partial M$.
Let $B^3$ denote the standard closed 3-ball
$\{x\in\bbr^3\mid |x|\le 1\}$.
An \emph{$n$-tangle} $T$ is a smooth pair
$(B^3, \lambda)$ where $\lambda$ is a proper embedding of
$n$ disjoint copies of the interval $I$ into $B^3$.
Throughout this paper, an embedding means either the map or the image.
Let $U_n$ denote a trivial $n$-tangle, \emph{i.e.}, $U_n$ consists of $n$
unlinked unknotted arcs.
For example, $U_1$ is the unknotted standard ball pair $(B^3,I)$.
For $n=2$, see Figure~\ref{fig:tangle}.

Concordances and invertible concordances between tangles
can be defined in a similar way as for links.
However, the boundary of the 3-ball $B^3$ is required to be fixed at each
stage of
concordance.  More precisely, let $I_1,\ldots,I_n$, denote $n$ disjoint
copies of the interval $I$.
Two $n$-tangles, $T_0=(B^3,\lambda_0)$ and $T_1=(B^3,\lambda_1)$,
are \emph{concordant} if there is a proper smooth embedding $\tau$ of
$(\cup_{i=1}^n I_i)\times I$ into $B^3\times I$, with
$\tau(\cup_{i=1}^n I_i\times\epsilon)=\lambda_\epsilon\ (\epsilon=0,1)$
and $\tau(\epsilon_i\times I)=\tau(\epsilon_i\times 0)\times I$
for each $i=1,\ldots, n$, and $\epsilon_i=0,1$ in $I_i$.
Let $(V;T_1,T_2)$ denote $(B^3\times I,\tau)$, the concordance between
$T_1$ and $T_2$.
If $(V;T_1,T_2)$ and $(V';T_2,T_3)$ are two concordances, we can then
paste $V'$ to $V$ along $T_2$ to get a concordance $(V\cup V';T_1,T_3)$.
A concordance $(V;T_1,T_2)$ is \emph{invertible at}
$T_2$ if there is a concordance $(V';T_2,T_1)$ such that $(V\cup
V';T_1,T_1)$ is diffeomorphic to $(T_1\times I;T_1,T_1)$ by a
diffeomorphism $\varphi$ with
$\varphi(\tau)=\lambda_1\times I$, where $\tau$ is the embedding of $n$
disjoint copies of $I\times I$ into $B^3\times I$ defining the concordance
$(V\cup V';T_1,T_1)$ and $\lambda_1$ is
the embedding of $n$ disjoint copies of $I$ into $B^3$ defining
the tangle $T_1$.

A knot is called \emph{doubly null concordant} if it is the slice of some
unknotted 2-sphere in $S^4$.  Two knots $K_1$ and $K_2$
are said to be \emph{doubly concordant} if $K_1\# J_1$ is isotopic
to $K_2\# J_2$ for some doubly null concordant knots $J_1$ and $J_2$.

The following theorem is due to Zeeman.
\begin{thm}  \cite{zee65}
\label{thm:zeeman}
Every 1-twist-spun knot is unknotted.
\end{thm}
Let $-K$ denote the knot obtained by taking the image of $K$, with
reversed orientation, under a reflection of $S^3$.
The following fact was first proved by Stallings and now follows readily
from \ref{thm:zeeman}.
(One cross-section of the 1-twist-spin of $K$ yields $K\# (-K)$. For
details, see \cite{sum71}.)
\begin{cor}
\label{cor:mirror}
$K\# (-K)$ is doubly null concordant for every knot $K$.
\end{cor}

\begin{cor}
If $K_1\# (-K_2)$ is doubly null concordant then $K_1$ and $K_2$ are
doubly concordant.
\end{cor}
\begin{proof}
Take $J_1 = K_2\# (-K_2)$ and $J_2= K_1\#(-K_2)$ in the definition of
double concordance.
\end{proof}

\begin{remark}
An easy exercise shows that knots $K_1$ and $K_2$ are concordant if and
only if $K_1\#(-K_2)$ is \emph{slice}, \emph{i.e.},
concordant to the unknot. This defines an equivalence relation.  However, a
definition of double concordance more along the lines of concordance
is as of yet
inaccessible.  The difficulty is that it is unknown whether the following
is true: If knots $K$ and $K\# J$ are doubly null concordant, then $J$ is
doubly null concordant. 
\end{remark}

There is a relation between invertible concordance and double concordance.

\begin{prop}
If $K_1$ is invertibly concordant to $K_2$ then $K_1\#(-K_2)$ is doubly
null concordant.
\end{prop}
\begin{proof}
There is a copy of $S^3\times I$ in $S^4$ intersecting the 1-twist-spin of
$K_1$ in $K_1\#(-K_1)\times I$. Since $K_2$ splits $K_1\times I$, there is
an invertible concordance from $K_1\#(-K_1)$ to $K_1\#(-K_2)$. Hence
$K_1\#(-K_1)\times I$ is split by $K_1\#(-K_2)$ and the result follows.
\end{proof}

\section{Invertible concordances and prime knots}
Kirby and Lickorish \cite{kl79} proved that any knot in
$S^3$ is concordant to a prime knot. Livingston \cite{liv81}
gave a different proof of this result using satellite knots.
In this section, we modify Livingston's approach to prove
Theorem~\ref{thm:main}.

Before proving this, we will set up some notation.
By a \emph{splitting}-$S^2$, $S$, for a knot $K$ (in $S^3$ or
$S^1\times D^2$) we denote an embedded 2-sphere, $S$, intersecting $K$
in exactly 2 points. A knot in either $S^3$ or $S^1\times D^2$ is
\emph{prime} if for every splitting-$S^2$, $S$, $S$ bounds some 3-ball,
$B$, with $(B,B\cap K)$ a trivial pair.  The \emph{winding number} of a
knot $K$ in $S^1\times D^2$ is that element $z$ of
$\bbz\cong H_1(S^1\times D^2;\bbz)$ with $z\ge 0$ and $K$ representing $z$. 
The \emph{wrapping number} of $K$ is the minimum number of intersections of
$K$ with a disk $D$ in $S^1\times D^2$ with $\partial D=$ meridian.
If $K_1$ is a knot in $S^1\times D^2$ and $K_2$ is a knot in $S^3$, the
$K_1$ \emph{satellite of} $K_2$ is the knot in $S^3$ formed by mapping
$S^1\times D^2$ into the regular neighborhood of $K_2$, $N(K_2)$,
and considering the image of $K_1$ under this map.  The only restriction
on the map of $S^1\times D^2$ into $N(K_2)$ is that it maps a meridian to
a meridian.
In what follows we will consider $S^1\times D^2$ embedded in $S^3$ in a
standard way. Hence any knot $K$ in $S^1\times D^2$ gives rise to a knot
$K^\ast$ in $S^3$.

The following theorem is due to Livingston.

\begin{thm}\cite{liv81}
\label{thm:chuck}
Let $K_1$ be a knot in $S^1\times D^2$ such that
$K_1^\ast$ is the unknot in $S^3$.
Then $K_1$ is prime in $S^1\times D^2$.
Moreover, if $K_1$ has wrapping number $>1$ and
$K_2$ is any nontrivial knot in $S^3$,
then the $K_1$ satellite of $K_2$ is prime in $S^3$.
\end{thm}

This theorem suggests that, to prove our main theorem~\ref{thm:main},
we only need to find a knot $K_1$ in $S^1\times D^2$ with $K_1^\ast$
the unknot in $S^3$ and an invertible concordance
between the core $C$ and the knot $K_1$ in $S^1\times D^2$.
To do this, we observe that there is an invertible concordance
between the tangles $U_2$ and $T$ in Figure~\ref{fig:tangle}.
We remark here that Ruberman in~\cite{rub90} has used the tangle $T$ to
prove that any closed orientable $3$-manifold is invertibly homology
cobordant to a hyperbolic $3$-manifold.

\begin{lem}
\label{lem:split}
The 2-tangle $T$ in Figure~\ref{fig:tangle}(b) splits $U_2\times I$.
\end{lem}

\begin{figure}
\input{fig1.pstex_t}
\caption{}
\label{fig:tangle}
\end{figure}

\begin{figure}
\input{fig2.pstex_t}
\caption{}
\label{fig:concordance}
\end{figure}

\begin{proof}
Let $I_1$ be a copy of the non-straight arc of $T$ in the 3-ball $B^3$
and let $J_1$ be a copy of the non-straight arc of $U_2$ in $B^3$
as shown in Figure~\ref{fig:tangle}(c).
The closed curve $J_1\cup I_1$ bounds an obvious punctured torus
$F$ that is the shaded region in Figure~\ref{fig:tangle}(c).
Consider $F$ as the plumbing of two $S^1\times I$. Let $c_i,\ i=1,2,$ be
the cores of the two $S^1\times I$ of $F$
and let $\bar{c}_i,\ i=1,2,$ be disjoint proper line segments
in~$F$ intersecting with $c_i$ exactly once, respectively. 
See Figure~\ref{fig:tangle}(c).

To construct an invertible concordance,
we will construct two concordances and then paste them together.
First, note that pinching $I_1$ along $\bar{c}_1$ transforms $T$
into the tangle $U_2$ with an unlinked unknotted circle inside which is
isotopic to the circle~$c_2$.
Now capping off this circle we have a concordance $(V_1';T,U_2)$.
The tangle $B^3\times \frac{1}{4}$ in Figure~\ref{fig:concordance}
represents a slice of this concordance before capping off the circle.
In the similar way, pinching $I_1$ along $\bar{c}_2$ and capping off the
unknot gives us another concordance $(V_2;T,U_2)$.
Let $(V_1;U_2,T)$ denote the concordance $(V'_1;T,U_2)$ with
reversed orientation.  We can then paste $V_1$ to $V_2$ along $T$ to
get a concordance $(V_1\cup V_2;U_2,U_2)$, which will be proved to be
isotopic to the product concordance $U_2\times I$.
A few cross-sections of concordance $V_1\cup V_2$
are drawn in Figure~\ref{fig:concordance}.

Let $\tau$ denote the embedding of two disjoint copies of $I\times I$ into
$V_1\cup V_2$ as in the definition of concordance
in Section~\ref{sec:definitions}.
It is obvious from Figure~\ref{fig:concordance} that
there is a 3-manifold $M$ (the union of shaded regions)
in $V_1\cup V_2$ bounded by $\tau$ and
$J_1\times I$, whose intersection with $U_2$
at each end of the concordance is the arc $J_1$ and
whose cross-section in the middle is the punctured torus~$F$.
This 3-manifold $M$ can be considered as the union of three submanifolds:
the product $F\times I$ and two 3-dimensional 2-handles $D^2\times I$.
One $D^2\times I$ is glued to $F\times I$ along a regular neighborhood of
$c_2$, which corresponds
to capping off the circle isotopic to $c_2$ as we constructed the
concordance $V'_1$.
The other $D^2\times I$ is glued along a regular neighborhood of $c_1$,
which corresponds to capping
off the circle isotopic to $c_1$ as we constructed the concordance $V_2$.
Since $F\times I$ is a 3-dimensional handlebody with 2 handles
with cores $c_1$ and $c_2$, $M$ is
the manifold that results by adding two 2-handles to a genus $2$ solid
handlebody along the cores of the 1-handles, in this case yielding $B^3$.
Moreover, $M$ does not intersect the other straight arc of $T$
at any stage.
Using this 3-ball $M$, we can isotop  $\tau$ to $J_1\times I$ in a regular
neighborhood of $M$ not disturbing the other arc and $\partial B^3$.
This completes the proof.
\end{proof}

\begin{prop}
\label{prop:split}
The knot $K_1$ in Figure~\ref{fig:satellite}(b) splits $C\times I$, where $C$ is
the core in $S^1\times D^2$.
\end{prop}

\begin{figure}
\input{fig3.pstex_t}

\bigskip

\input{fig4.pstex_t}
\caption{}
\label{fig:satellite}
\end{figure}

\begin{proof}
Consider $S^1\times D^2$ as the complement of the unknot $m$ in $S^3$.
The knot $K_1$ in Figure~\ref{fig:satellite}(b) is isotopic
to $K_1$ in Figure~\ref{fig:satellite}(a). It is
obvious from Figure~\ref{fig:satellite}(a) that $K_1\cup m$
is the link in $S^3$ formed by replacing a trivial 2-tangle in Hopf link
with $T$ (dotted circle in Figure~\ref{fig:satellite}(a)).
The proposition now follows from Lemma~\ref{lem:split}.
\end{proof}

Now we are ready to prove our main theorem~\ref{thm:main}.

\begin{proof}[Proof of Theorem~\ref{thm:main}]
Let $K$ be a knot in $S^3$.
If $K$ is trivial it is prime itself. Suppose now that $K$ is nontrivial.
Let $K'$ be $K_1$ satellite of $K$ where $K_1$ is the knot in $S^1\times D^2$
in Figure~\ref{fig:satellite}(b). By Proposition~\ref{prop:split}, $K'$ splits
$K\times I$. We now only need to show $K'$ is prime.  Since $K_1^\ast$ is the
unknot in $S^3$, $K_1$ is prime by Theorem~\ref{thm:chuck} and to complete
proof it remains to
show its wrapping number $> 1$.  Its winding number is 1, hence
its wrapping number is at least one.  It is easy to see that 
the only prime knot in $S^1\times D^2$ with
wrapping number~1 is the core.  So, if $K_1$ had wrapping number~1, then
it is isotopic to the core of $S^1\times D^2$.
The $-1$ surgery on the meridian curve $m$ in $S^3$ should make
$K_1^\ast$ unchanged, \emph{i.e.}, unknotted.
However, the knot in Figure~\ref{fig:satellite}(e),
the result of $K_1^\ast$ after $-1$ surgery along $m$,
is $9_{46}$ and hence knotted. Therefore the wrapping number is $>1$.
\end{proof}

\begin{cor}
Any knot is doubly concordant to a prime knot.
\end{cor}

\begin{remark}
The $K_1$ satellite of $K$ has the same Alexander polynomial as that
of~$K$. Seifert~\cite{sei50}
proved that the Alexander polynomial of the $K_1$
satellite of $K$ is $\Delta_{K_1^\ast}(t) \Delta_K(t^w)$
if $w$ is the winding number of $K_1$
in $S^1\times D^2$.  In our case, $w$ is $1$ and $K_1^\ast$ is the unknot.
\end{remark}

In \cite{liv81}, Livingston also proved that every 3-manifold is homology
cobordant to an irreducible 3-manifold.
Two 3-manifolds, $M_1$ and $M_2$, are \emph{homology cobordant} if there
is a 4-manifold $W$, with $\partial W=M_1\cup M_2$ and the map of
$H_\ast(M_i;\bbz)\to H_\ast(W;\bbz)$ an isomorphism.  Invertible homology
cobordisms can be defined in the same way as in the knot concordance case.
A 3-manifold $M$ is \emph{irreducible} if every embedded
$S^2$ in $M$ bounds an embedded~$B^3$. 

\begin{remark}
In spirit of \cite{liv81}, we have a simple proof
that every 3-manifold is invertibly
homology cobordant to an irreducible 3-manifold.
To prove this, we only need to slightly modify the proof of
Theorem~3.2 in \cite{liv81} by using $K_1$ in
Figure~\ref{fig:satellite}(b).
The $-1$ surgery on $K_1$ makes the meridian $m$ the knot $9_{46}$.
\end{remark}

This remark is also a corollary of Ruberman's Theorem~2.6 in
\cite{rub90} that reads:
for every closed orientable 3-manifold $N$, there is a
hyperbolic 3-manifold $M$, and an invertible homology cobordism from $M$
to $N$. The remark follows since a hyperbolic 3-manifold is
irreducible.

\end{document}

%% file: fig1.pstex_t
\begin{picture}(0,0)%
\epsfig{file=fig1.pstex}%
\end{picture}%
\setlength{\unitlength}{908sp}%
\begin{picture}(23291,7364)(10186,-8161)
\put(27301,-5611){\makebox(0,0)[b]{\tiny$F$}}
\put(28951,-8161){\makebox(0,0)[b]{\small (c)}}
\put(21076,-8161){\makebox(0,0)[b]{\small (b) $T$}}
\put(13351,-8161){\makebox(0,0)[b]{\small (a) $U_2$}}
\put(27676,-1711){\makebox(0,0)[b]{\tiny$\bar{c}_1$}}
\put(29026,-1711){\makebox(0,0)[b]{\tiny$\bar{c}_2$}}
\put(26326,-3211){\makebox(0,0)[b]{\tiny$c_1$}}
\put(31576,-3136){\makebox(0,0)[b]{\tiny$c_2$}}
\put(32701,-3961){\makebox(0,0)[b]{\tiny$I_1$}}
\put(29176,-6061){\makebox(0,0)[b]{\tiny $J_1$}}
\end{picture}

%% file: fig2.pstex_t
\begin{picture}(0,0)%
\epsfig{file=fig2.pstex}%
\end{picture}%
\setlength{\unitlength}{908sp}%
\begin{picture}(23130,7322)(1636,-8119)
\put(4801,-8011){\makebox(0,0)[b]{\scriptsize$B^3\times\frac{1}{4}$}}
\put(13201,-8011){\makebox(0,0)[b]{\scriptsize$B^3\times\frac{1}{2}$}}
\put(21601,-8011){\makebox(0,0)[b]{\scriptsize$B^3\times\frac{3}{4}$}}
\end{picture}

%% file: fig3.pstex_t
\begin{picture}(0,0)%
\epsfig{file=fig3.pstex}%
\end{picture}%
\setlength{\unitlength}{1184sp}%
\begin{picture}(17970,6435)(5931,-11461)
\put(6901,-6211){\makebox(0,0)[b]{\small$K_1$}}
\put(11401,-5911){\makebox(0,0)[b]{\small$m$}}
\put(19801,-11461){\makebox(0,0)[b]{\small(b)}}
\put(11026,-11461){\makebox(0,0)[b]{\small(a)}}
\put(21901,-5386){\makebox(0,0)[b]{\small$K_1$}}
\put(23101,-9061){\makebox(0,0)[b]{\small$m$}}
\end{picture}

%% file: fig4.pstex_t
\begin{picture}(0,0)%
\epsfig{file=fig4.pstex}%
\end{picture}%
\setlength{\unitlength}{1302sp}%
\begin{picture}(14913,5034)(2539,-5611)
\put(5401,-1411){\makebox(0,0)[b]{\small$K_1$}}
\put(9826,-1411){\makebox(0,0)[b]{\small$m$}}
\put(14476,-5611){\makebox(0,0)[b]{\small(e)}}
\put(8251,-5611){\makebox(0,0)[b]{\small(d)}}
\put(3901,-5611){\makebox(0,0)[b]{\small(c)}}
\put(6151,-2836){\makebox(0,0)[b]{\small$m$}}
\put(10576,-2836){\makebox(0,0)[b]{\small$K_1$}}
\put(14476,-2911){\makebox(0,0)[b]{$=$}}
\put(16276,-5161){\makebox(0,0)[b]{\small$9_{46}$}}
\end{picture}